\documentclass[11pt,reqno]{amsart}

\usepackage{amsmath}
\usepackage{amsfonts}
\usepackage{amssymb}
\usepackage{amsthm,color}
\usepackage{amsmath,amsthm,amssymb}
\usepackage{graphicx,mathtools}
\usepackage{cite}
\newtheorem{Theorem}{Theorem}[section]
\newtheorem{Definition}[Theorem]{Definition}
\newtheorem{Proposition}[Theorem]{Proposition}
\newtheorem{Lemma}[Theorem]{Lemma}
\newtheorem{Corollary}[Theorem]{Corollary}

\def\a{\alpha}
\def\d{\delta}

\def\({\left(}
\def\){\right)}

\def\t{\theta}
\def\g{\gamma}

\newcommand{\E}{\mathbf E}

\begin{document}

\title[An inverse potential problem]{An inverse potential problem for the stochastic heat equation with space-time noise}

\author{Peijun Li}
\address{LSEC, ICMSEC, Academy of Mathematics and Systems Science, Chinese Academy of Sciences, Beijing 100190, China}
\email{lipeijun@lsec.cc.ac.cn}  

\author{Xiangchan Zhu}
\address{Academy of Mathematics and Systems Science, Chinese Academy of Sciences, Beijing 100190,
China} 
\email{zhuxiangchan@126.com}

\author{Yichun Zhu}
\address{Academy of Mathematics and Systems Science, Chinese Academy of Sciences, Beijing 100190,
China} 
\email{steven00931002@hotmail.com}

\thanks{The work of PL is supported by the National Key R\&D Program of China (2024YFA1012300).}



\begin{abstract}
This paper investigates an inverse potential problem for the stochastic heat equation driven by space-time Gaussian noise, which is spatially colored and temporally white. The objective is to determine the covariance operator of the random potential. We establish that the covariance operator can be uniquely identified from the correlation of the mild solution to the stochastic heat equation at a final time, where the initial conditions are specified by a complete orthonormal basis. The analysis relies on characterizing a tensor product structure inherent in the problem and utilizing the monotonicity properties of the operators associated with the system.
\end{abstract}

\maketitle

\section{Introduction}

In this paper, we consider the stochastic heat equation driven by space-time Gaussian noise, formulated as follows:
\begin{equation}\label{eq:SHE}
\begin{cases}
{\partial_t} u(t,x)= \Delta u(t,x) + u(t,x)\xi(t,x),\quad &(t,x)\in (0, T]\times D,\\
u(0,x)=u_0(x),\quad & x\in D,\\
u(t, x)=0,\quad & (t, x)\in (0, T]\times \partial D,
\end{cases}
\end{equation}
where $D \subseteq \mathbb{R}^d$ is a bounded domain with a Lipschitz boundary $\partial D$,  and $\xi(t,x)$ represents Gaussian noise, with its regularity and properties to be specified later. The initial condition $u_0(x)$ is assumed to satisfy $u_0\in L^2(D)$, and the solution is subject to a homogeneous Dirichlet boundary condition. The potential $\xi(t,x)$ introduces stochastic forcing, which models random fluctuations in the system. This work addresses the inverse problem of identifying the random potential $\xi$ based on the statistical properties of the solution $u$. 

The inverse potential problem for the heat equation, and more generally for the diffusion equation, including the 
time-fractional diffusion equation, has been extensively studied in deterministic settings \cite{CK-SIMA01, CNYY-IP09, JZ-IP21, KR-IP19, KOSY-JDE18}. We refer to the monograph \cite{I-98} for a comprehensive overview of inverse problems associated with the diffusion equation. In contrast, there has been relatively little research addressing the inverse problems for the stochastic diffusion equation. To the best of our knowledge, the only available work is presented in \cite{FLW-IPI24}. Specifically, the authors consider the stochastic diffusion equation driven by a multiplicative white noise, expressed as $u(t,x)\xi(t,x)=u(t,x)q(t)\dot{W}(t)$, where $W(t)$ denotes the one-dimensional standard Brownian motion and the time-dependent function $q(t)$ represents the strength of the randomness. The inverse problem is to determine the potential function $q$ from measurements $u$ taken at an interior observation point within the spatial domain over a time interval. Based on the analytical solution to the stochastic diffusion equation, an explicit formula is derived to reconstruct $q^2$ using the expectation of the measurements over the sample space. 

In this work, we consider a more general class of Gaussian noise driven by a cylindrical Wiener process, which is spatially colored and temporally white. The approach differs significantly from the technique used in \cite{FLW-IPI24}. For the direct problem, it follows from \cite{LR} that the stochastic heat equation \eqref{eq:SHE} has a unique mild solution (cf. Theorem \ref{thm:wellposed}). For the inverse problem, we demonstrate that the covariance operator associated with the random potential can be uniquely determined from the correlation of the mild solution to the stochastic heat equation \eqref{eq:SHE} at a final time, where the initial conditions given by a complete orthonormal system in $L^2(D)$. The key components of the proof involve identifying a tensor product structure and utilizing the $C_0$ semigroup generated by a monotone operator.

It is worth noting that there has been considerable recent attention to the inverse random potential problem for time-harmonic wave equations \cite{LPS-CMP08, LLW-SIMA22, LW-SIMA24}. In this context, the potential is modeled as a microlocally isotropic Gaussian rough field. The inverse problem focuses on estimating the correlation strength of the random potential by studying the second moment of the scattering data in the high frequency limit, averaged across a frequency range. However, research on the inverse potential problem for stochastic wave equations in the time domain remains limited. This gap presents an opportunity for future exploration, building upon the approach developed in this work.

The structure of the paper is as follows. Section \ref{sec:introduction} outlines the problem formulation and summarizes the main results of the study. Section \ref{sec:functional} introduces the necessary functional preliminaries. Section \ref{sec:inverse} establishes the uniqueness of the solution to the inverse problem. The paper concludes with general remarks in Section \ref{sec:conclusion}.

\section{Problem formulation and the main result}\label{sec:introduction}

Let $H$ denote the standard $L^2$ space over $D$, equipped with the norm $\|\cdot\|$ and inner product $\langle \cdot, \cdot\rangle$. Let $A$ be the Laplacian on $D$ with the Dirichlet boundary condition, and let $\{S(t)\}_{t \geq 0}$ denote the analytic semigroup generated by $A$. Since $D$ is a bounded domain, there exists a complete orthonormal system $\{e_k\}_{k \in \mathbb{N}}$ in $H$, along with a strictly increasing sequence of positive numbers $\{\a_k\}_{k \in \mathbb{N}}$, such that 
\begin{equation}\label{Aek}
-Ae_k = \a_k e_k,\quad k \in \mathbb{N}.
\end{equation}

For $s\geq 0$, $p>1$, define
\begin{equation*}
W^{s,p}:= \{ f \in L^p(D): \|(\mathbb{I}-\Delta)^{s/2}f\|_{L^p}< \infty \},
\end{equation*}
where $\mathbb{I}$ is the identity operator. Let $H^{\gamma}$ denote the completion of $C^{\infty}_c(D)$ in $W^{\gamma,2}$, and let $H^{-\gamma}$ be the dual space of $H^{\gamma}$. It is known that $H^{2\gamma}$ and $D((-A)^{\gamma})$ have equivalent norms. Let $U_1$ and $U_2$ be two Hilbert spaces. We denote by $\mathcal{L}_2(U_1;U_2)$ the Hilbert--Schmidt operator from $U_1$ to $U_2$. 

In the following, we fix a complete probability space $(\Omega, \mathcal{F}, \mathbf{P})$. Define an operator $Q$ on $H$ by 
\[
Qe_k= \lambda_k e_k,\quad k \in \mathbb{N},
\]
where $\{\lambda_k\}_{k\in\mathbb N}$ are the eigenvalues of $Q$. The cylindrical Wiener process on $H$ is defined as
\[
W_t:= \sum_{k \in \mathbb{N}} e_k \beta_k(t),
\]
where $\{\beta_k\}_{k\in\mathbb N}$ is a collection of mutually independent one-dimensional standard Brownian motions on $(\Omega, \mathcal{F}, \mathbf{P})$. We define $\{\mathcal{F}_t\}_{t \geq 0}$ as the normal filtration generated by $W_t$. The random potential $\xi(t,x)$, which is spatially colored and temporally white, is given by
\begin{equation*}
\xi(t,x): =Q \partial_tW_t= \sum_{k \in \mathbb{N}} \lambda_k e_k(x) d\beta_k(t),
\end{equation*}
where the eigenvalues $\lambda_k$ determine the intensity of the corresponding spatial modes $e_k$ in the random potential $\xi(t, x)$. These eigenvalues describe how each spatial mode contributes to the overall behavior of the random potential. Consequently, the random potential $\xi(t, x)$ is fully characterized by the operator $Q$, which is determined by the covariance eigenvalues $\lambda_k$ for $k\in\mathbb N$.

With these notations, the stochastic heat equation \eqref{eq:SHE} can be expressed in its abstract form as
\begin{equation*}%\label{eq:SHE:0}
d u(t)=A u(t,x) + G(u(t)) Q dW_t,\quad u(0)=u_0, 
\end{equation*}
where $G: L^{d/(2\g)} \to \mathcal{L}(H^{2\g}, H)$, $\g \in [0, d/4)$, is the Nemytskii operator defined by 
\begin{equation*}
G(u)f(x):= u(x)f(x),\quad  x \in D.
\end{equation*}
By Sobolev embedding, we obtain
\begin{equation}\label{eq:G}
\|G(u)f\|\leq \|f\|_{L^{2d/(d-4\g)}} \|u\|_{L^{d/(2\g)}} \leq \|f\|_{H^{2\g}} \|u\|_{L^{d/(2\g)}},\quad \g \in [0, d/4),
\end{equation}
where the first inequality follows from H\"{o}lder's inequality. 

For $\gamma>0$ and the operator $Q$, define 
\begin{equation*}
  \Lambda_{\gamma,Q} : =\sum_{k \in \mathbb{N}}{\lambda}_k^2 \|e_k\|_{L^{d/2\gamma}}^2.
\end{equation*}
This quantity plays a critical role in ensuring the well-posedness of \eqref{eq:SHE} and the uniqueness of the solution to the inverse problem. We refer to Section \ref{r:data} for a more detailed discussion on the condition $\Lambda_{\gamma,Q} < \infty$.

\subsection{Main results}

We begin with the definition of a solution to the stochastic heat equation \eqref{eq:SHE}. 

\begin{Definition}
Let $T>0$ be an arbitrary finite time. A continuous $H$-valued process $\{u(t)\}_{t \in [0,T]}$ is said to be a solution to the stochastic heat equation \eqref{eq:SHE} if $u \in L^2([0,T]\times \Omega, dt \times d\mathbf{P};H^1)$ and satisfies, for any $t \in [0,T]$, the integral equation
\begin{equation*}%\label{eq:mild}
 u(t)= u_0 +\int_0^t A u(s) ds+ \int_0^t  G(u(s)) QdW_s,\quad {\mathbf P}\text{-a.s.}
 \end{equation*}
\end{Definition}

For any $u,v \in H^1$ and $\g \in [0, d/4)$, by \eqref{eq:G}, we have
\begin{align*}
&\langle A(u-v),u-v \rangle + \|(G(u)-G(v))Q\|_{\mathcal{L}_2^0(H;H)}^2 \\
&\leq -\|u-v\|_{H^1}^2 + \sum_{k \in \mathbb{N}} \lambda_{k}^2\ \|(u-v)e_k\|_{H}^2 \\
&\leq -\|u-v\|_{H^1}^2 + \sum_{k \in \mathbb{N}} \lambda_{k}^2\ \|e_k\|_{L^{d/(2\g)}}^2 \|u-v\|_{H^{2\g}}^2.
\end{align*}
When $\g \in [0, \frac{1}{2}\wedge \frac{d}{4})$, by interpolation and Young's inequality, there exists a positive constant $C$ such that
\begin{equation*}
\langle A(u-v),u-v \rangle + \|(G(u)-G(v))Q\|_{\mathcal{L}^2_0(H;H)}^{2}\leq- \frac{1}{2}\|u-v\|_{H^1}^2 + C \Lambda_{\g,Q}\  \|u-v\|^2.
\end{equation*}
Moreover, the mapping $t \to \langle A(u+tv), w \rangle$ is continuous for all $u,v,w \in H^1$. Therefore, applying the results in \cite[Theorem 4.2.5]{LR}, we obtain the following theorem.

\begin{Theorem}\label{thm:wellposed}
 Let $Q$ be such that there exists $\g  \in [0, \frac{1}{2}\wedge \frac{d}{4})$ satisfying $ \Lambda_{\g,Q}< \infty$. Then, there exists a unique mild solution $u$ to the stochastic heat equation \eqref{eq:SHE} such that
\begin{equation}\label{eq:u:reg}
\int_0^T \E[\|u(t)\|_{H^1}^2] dt + \E[\sup_{t \in [0,T]} \|u(t)\|^2] <\infty.
\end{equation}
\end{Theorem}

To illustrate our main result, we introduce the symmetric Hilbert space $U$, which is a subspace of $L^2(D\times D)$ defined by 
\begin{equation*}
U:= \{f : f(x,y)=f(y,x),\, \|f\|_U^2:=\int_{D \times D} |f(x,y)|^2\ dx dy< \infty\}, 
\end{equation*}
 with the inner product $\langle \cdot, \cdot\rangle_U$ induced by the norm. Let  $A_0$ denote the Laplacian on $D\times D$ with the Dirichlet boundary condition in $U$. The domain of $A_0$ is given by $D(A_0):= \{x \in U: A_0x \in U\}$. For any $\g \in \mathbb{R}$, we define the space $U^{\g}$ as $D((-A_0)^{\g/2})\cap U$.
  
Henceforth, we denote $f\otimes g$ by $f(x) g(y)$ for all $f,g \in H$. Using these notations, we formally define $\theta(t,u_0)$ as
 \[
 \t(t,u_0):= \E[u(t,u_0) \otimes u(t,u_0)],
 \] 
where, as stated in Theorem \ref{thm:wellposed}, $u \in L^{\infty}(0,T;H)$. Moreover, we will show in Lemma \ref{l:u:regularity} that $\theta$ belongs to  $L^{\infty}(0,T;U) \cap L^2(0,T;U^1)$ for all $u_0 \in H$.

The main result of this paper is presented in the following theorem, which addresses the uniqueness of the inverse problem. 

\begin{Theorem}\label{thm:main}
Assume that there exists some $\gamma \in [0, \frac{1}{2}\wedge \frac{d}{4})$ such that $\Lambda_{\g,Q}<\infty$. For any $i, j \in \mathbb{N}$, we define $\theta^{i,j}(t)$ as
\begin{equation*}
\t^{i,j}(t)= \theta(t, e_i+e_j)-\theta(t, e_i-e_j). 
\end{equation*}
Then, the covariance operator $Q$  can be uniquely determined by $\{\t^{i,j}(t_0)\}_{i,j\in\mathbb N}$ at some given fixed time $t_0>0$. 
\end{Theorem}

\subsection{Remarks on random potential}\label{r:data}

As demonstrated in the proof, we assume throughout the paper that $\Lambda_{\gamma, Q} < \infty$ for some $\gamma \in [0,\frac{1}{2}\wedge \frac{d}{4})$. This condition characterizes the regularity of the random potential $\xi(t, x)$. The uniqueness of the inverse problem holds when $\xi(t, x)$ has sufficient regularity.

When $\g_1>\g_2$, by H\"{o}lder's inequality, there exists a positive constant $C_{D,\g_1,\g_2}$, independent of $k$, such that
\[
\|e_k\|_{L^{d/(2\g_1)}} \leq C_{D,\g_1,\g_2} \|e_k\|_{L^{d/(2\g_2)}}.
\]
Hence, we have $\Lambda_{\g_1,Q} \leq C_{D,\g_1,\g_2} \Lambda_{\g_2,Q}$ for $\g_1\geq \g_2$, and our condition is more general than 
\[
\sum_{k \in \mathbb{N}} \lambda_k^2\|e_k\|_{L^{\infty}}^{2}  < \infty.
\]

\subsection{Remarks on data}

Theorem \ref{thm:main} states that the covariance of the random potential can be uniquely determined from the covariance of the solution for different choices of the initial condition. We note that it is unlikely that fewer data would suffice to reconstruct $Q$.

Let $\gamma \in[0,\frac{1}{2}\wedge \frac{d}{4})$, and let $H_Q:U^{4\gamma} \to U$ be an operator defined by 
\begin{equation}\label{eq:HQ}
H_Qf:= \left(\sum_{k \in \mathbb{N}} \lambda_k^2\ e_k \otimes e_k \right) f,
\end{equation}
which will be specified in Lemma \ref{l:HQ}. It is clear that the information about $Q$ is fully encapsulated in the operator $H_Q$, which can be deduced from the structure of the semigroup $e^{(H_Q+A_0)t}$ and the data $\theta^{i,j}_{m,n}$. The problem becomes significantly simpler when $H_Q$ and $A_0$ commute, as this allows the factorization
\[
e^{(H_Q+A_0)t}= e^{H_Qt}e^{A_0t}. 
\]
This factorization enables the decomposition of the dynamics into two separate linear systems:
\[
d \theta(t)= A_0 \theta(t),\quad d \theta(t)= H_Q \theta(t).
\]
Such a decomposition reduces the amount of data needed to recover the covariance matrix $Q$. 

In the general case where $H_Q$ and $A_0$ do not commute, multiple data sets generated from different initial conditions are required, as indicated by Theorem \ref{thm:main}. According to Trotter's theorem \cite[Corollary 5.5]{P}, we have
\[
e^{(H_Q+A_0)t}x= \lim_{n \to \infty} \big(e^{H_Q t/n} e^{A_0 t/n}\big)^nx,\ \ \ x \in H.
\]
This product formula highlights the intertwining of the two operators, suggesting that additional effort is necessary to disentangle the contributions of $H_Q$ and $A_0$. Furthermore, when determining a Gaussian measure on the Hilbert space, it is essential to know every entry of the covariance operator. Finally, due to the inherent symmetry in the choice of $e_i+e_j$, it is unlikely that fewer data would suffice to reconstruct $Q$.

\subsection{Idea of the proof}
The main idea of the proof is to show that $\theta$ satisfies
\begin{equation*}%\label{eq:strong:0}
\frac{d}{dt} \t(t)= (H_Q+A_0) \t(t).
\end{equation*}
 By solving this equation, we obtain an explicit formula for $\theta(t)$:
\[
\t(t)= e^{(H_Q+A_0)t}\t(0).
\]
By  Lemma \ref{l:C0}, $A_0+H_Q$ generates a $C_0$ semigroup. The key step in the proof is that we can recover the operator $e^{(H_Q+A_0)t_0}$ from the given data. More precisely, by the analytic property of the semigroup and the compactness property of self-adjoint operator, we can derive all the discrete spectrum of $\{e^{(H_Q+A_0)t}\}_{t>0}$ at any time $t>0$. Therefore, the operator $H_Q$ can be uniquely determined from the given data.

\section{Functional preliminaries}\label{sec:functional}

This section is dedicated to studying the properties of $A_0+H_Q$. In the following, we will prove that $A_0+H_Q$ is a self-adjoint operator that generates a $C_0$ semigroup $\{e^{(A_0+H_Q)t}\}_{t \geq 0}$. We begin by examining some properties of the heat semigroup generated by $A_0$.

 Let $\{S_0(t)\}_{t >0}$ denote the analytic semigroup generated by the Laplacian $A_0$. It is known that $\{S_0(t)\}_{t \geq 0}$ forms an analytic $C_0$-semigroup on $U$, and there exists a positive constant $C_{\gamma_1,\gamma_2}$ such that the following smoothing property of the heat semigroup holds:
 \begin{equation}\label{eq:Schauder}
 \|S_0(t)x\|_{U^{\g_1}} \leq C_{\gamma_1,\gamma_2}\ t^{-(\gamma_1-\gamma_2)/2}\|x\|_{U^{\gamma_2}},\quad\forall \gamma_1 > \gamma_2.
 \end{equation}

The following important characterization of $U^{\g}$ is provided in \cite[p.100, ex.3]{KF}.

\begin{Proposition}\label{p:HH}
The set $\{e_i(x)e_j(y)\}_{i,j \in \mathbb{N}}$ forms a complete orthonormal system in $U$. Moreover, for any $\g \in \mathbb{R}$, the following characterization holds:
\begin{equation*}
U^{\g}:= \{f=\sum_{i,j \in \mathbb{N}} f_{i,j} e_i(x)  e_j(y):  f_{i,j}=f_{j,i}, \, \|f\|_{U^{\g}}^2 = \sum_{i,j \in \mathbb{N}} \a_i^{\g} \a_j^{\g} |f_{i,j}|^2 < \infty  \},
\end{equation*}
where $\alpha_i$, for $i\in\mathbb N$, is defined in \eqref{Aek}. 
\end{Proposition}

The following lemma provides a definition for the operator $H_Q$ introduced in \eqref{eq:HQ}.

\begin{Lemma}\label{l:HQ}
Let $Q$ and $\gamma_1, \g_2 \in [0, d)$ be such that $\Lambda_{(\g_1+\g_2)/4,Q}<\infty$. We  define $H_{Q}^{}: U^{\g_1}\to U^{-\g_2}$ by 
\begin{equation*}
H_Q^{}f:= \left(\sum_{k \in \mathbb{N}} \lambda_k^2\ e_k \otimes e_k \right) f,\quad \forall f \in U^{\gamma_1},
\end{equation*}
 and it follows that
\begin{equation*}
\|H_Q^{} \|_{\mathcal{L}(U^{\g_1};U^{-\g_2})} \leq \Lambda_{(\g_1+\g_2)/4,Q}.
\end{equation*}
\end{Lemma}

\begin{proof}
First, we observe that $H_Qf \in U$ for all  $f \in C^{\infty}_{c}(D \times D)\cap U$. It follows from Sobolev embedding and duality that 
\begin{equation*}
\|H_Qf\|_{U^{-\g_2}} \leq \|H_Qf\|_{L^{2d/(d+\g_2)}} \leq \sum_{k \in \mathbb{N}} \lambda_k^2 \|(e_k \otimes e_k) f\|_{L^{2d/(d+\g_2)}}.
\end{equation*}
Using H\"{o}lder's inequality, we have 
\begin{align*}
\|H_Qf\|_{U^{-\g_2}} &\leq \sum_{k \in \mathbb{N}} \lambda_{k}^2 \(\int_{D \times D}\big|  e_k(x) e_k(y)\big|^{2d/(\g_1+\g_2)} dx\ dy\)^{(\g_1+\g_2)/2d}\\
&\qquad \times \(\int_{D \times D}|  f(x,y)|^{2d/(d-\g_1)} dx\ dy\)^{(d-\g_1)/(2d)}.
\end{align*}
Applying Sobolev embedding again, we obtain $\|f\|_{L^{2d/(d-\g_1)}(D \times D)} \lesssim \|f\|_{U^{\g_1}}$, and thus
\begin{equation*}
\|H_Qf\|_{U^{-\g_2}} \leq \sum_{k \in \mathbb{N}} \lambda_{k}^2 \|e_k\|_{L^{2d/(\g_1+\g_2)}(D)}^2 \|f\|_{U^{\g_1}}=\Lambda_{(\g_1+\g_2)/4,Q}  \|f\|_{U^{\g_1}},
\end{equation*}
which completes the proof.
\end{proof}

\begin{Lemma}\label{l:selfad}
Let $Q$ be such that $\Lambda_{\g,Q}<\infty$ for some $\gamma \in [0, 1/2)$. Then $A_0+H_Q:D(A_0) \to U$ is self-adjoint and has a compact resolvent.
\end{Lemma}

\begin{proof}
It holds that for any $\varepsilon>0$, there exists a positive constant $C_{\varepsilon}$ such that
\[
\begin{aligned}\label{eq:res}
 \langle  H_Q f, f \rangle_U &= \sum_{k \in \mathbb{N}} \lambda_k^2 \int_{D \times D} e_k(x) e_k(y) |f(x,y)|^2 dx dy\notag\\
 &\leq \sum_{k \in \mathbb{N}} \lambda_k^2  \|e_k\|_{L^{d/2\g}(D)}^2\ \|f\|_{L^{2d/(d-2\g)}(D \times D)}^2\notag \\
 &\leq \sum_{k \in \mathbb{N}} \lambda_k^2  \|e_k\|_{L^{d/2\g}(D)}^2\ \|f\|_{U}^{1-2\gamma} \|f\|_{U^1}^{2\gamma}\notag \\
 &\leq   \varepsilon  \| f  \|_{U^1}^2 + C_{\varepsilon}  \Lambda_{\g,Q}  \|f\|_{U}^{2} .
\end{aligned}
\]
Here, in the second inequality, we apply H\"{o}lder's inequality; in the third inequality, we use the Sobolev inequality and interpolation; and in the final step, we apply Young's inequality.

 By \cite[Theorem XIII.68]{RS}, it follows that $A_0+H_Q$ has a compact resolvent.  By Lemma \ref{l:HQ}, $H_Q:U^{2\gamma} \to U^{-2\gamma}$ is bounded operator with norm $\Lambda_{\g,Q}<\infty$. Since $U^{2\gamma}$ is dense in $U$ and for all $f,g \in U^{2\gamma}$, there holds
\[
\langle H_Q f,g \rangle_U = \langle f, H_Qg \rangle_U,
\]
we obtain that $H_Q$ is a symmetric operator. By \eqref{eq:res} again and Kato--Rellich Theorem, $A_0+H_Q$ is self-adjoint with domain $D(A_0)$.
\end{proof}

\begin{Lemma}\label{l:C0}
Let $Q$ be such that $\Lambda_{\g,Q}<\infty$ for some $\gamma \in [0, 1/2)$. The operator $A_0+H_Q$ generates a $C_0$-semigroup $\{e^{(A_0+H_Q)t}\}_{t \geq 0}$ on $U$. 
\end{Lemma}

\begin{proof}
 Using \eqref{eq:res} yields 
\begin{equation*}%\label{eq:dissipative}
 \langle  (A_0+H_Q) f, f \rangle_U \leq  -\frac{1}{2}\langle \nabla f , \nabla f \rangle_U + C_{1/2}\Lambda_{\g,Q}  \|f\|_{U}^{2},\quad\forall f \in D(A_0),
\end{equation*}
which implies that $A_0+H_Q$ is quasi-dissipative. The result can then be deduced in the standard way. For the reader's convenience, we provide the proof below.

Let $\omega>C_{1/2}\Lambda_{\g,Q}$. It suffices to prove that the range of $\omega\mathbb{I}-A_0-H_Q$ is dense in $U$. If this does not hold, then there exists $f\neq 0$ in $U$ such that
\[
\langle (\omega \mathbb{I} -A_0-H_Q)u, f \rangle_U = 0,\quad \forall u \in D(A_0+H_Q).
\]
This implies that $f \in D(A_0)$. Since $D(A_0)$ is dense in $U$ and by Lemma \ref{l:selfad}, we have $(\omega \mathbb{I} - A_0-H_Q)f=0$, which implies
\[
0 = \langle (\omega \mathbb{I} - A_0-H_Q)f ,f \rangle \geq \frac{1}{2}\langle \nabla f , \nabla f \rangle_U + (\omega-C_{1/2}\Lambda_{\g,Q})  \|f\|_{U}^{2} \geq 0.
\]
This leads to  $f=0$, which is a contradiction. Hence, the range of $\omega\mathbb{I}-A_0-H_Q$ is dense in $U$, and the result follows by the Lumer--Phillips theorem \cite[Chapter 1.4]{P}. 
\end{proof}

\begin{Lemma}\label{l:analytic}
Let $Q$ be such that $\Lambda_{\g,Q}<\infty$ for some $\gamma \in [0, 1/2)$. Then the semigroup $\{e^{(A_0+H_Q)t}\}_{t \geq 0}$ is analytic. Moreover, $e^{(A_0+H_Q)t}$ is self-adjoint and compact on $U$.
\end{Lemma}

\begin{proof}
By Lemma \ref{l:C0}, $A_0 + H_Q$ is quasi-dissipative, and there exists a constant $c$ such that $[c, +\infty)$ belongs to the resolvent set of $A_0 + H_Q$. Since $A_0+H_Q$ is self-adjoint, it follows from \cite[Corollary 10.6]{P} that $e^{(A_0+H_Q)t}$ is also self-adjoint.

Moreover, by \cite[Chapter VII.2]{RS} that there exists a unique spectral family $\{E_{\lambda}\}_{\lambda\in{\mathbb{R}}}$ such that
\begin{equation*}
(A_0+H_Q) e^{(A_0+H_Q)t}= \int_{-\infty}^{c} \lambda e^{\lambda t}dE_{\lambda},
\end{equation*}
and
\begin{equation*}
\|(A_0+H_Q)e^{(A_0+H_Q)t} \|_{\mathcal{L}(U)} = \sup_{\lambda \in (-\infty,c]} |e^{\lambda t}\lambda|.
\end{equation*}
Therefore, for any $T>0$, there exists a positive constant $C_T$ such that
\begin{equation*}
\|(A_0+H_Q)e^{(A_0+H_Q)t} \|_{\mathcal{L}(U)}\leq \frac{C_T}{t},\ \ t \in[0,T].
\end{equation*}
By \cite[Theorem 5.2]{P},  $\{e^{(A_0+H_Q)t}\}_{t \geq 0}$ is an analytic semigroup on $U$. 

To prove that $e^{(A_0+H_Q)t}$ is compact for $t >0$, we observe the following identity:
\begin{equation*}
e^{(A_0+H_Q)t}= (\lambda- (A_0+H_Q))e^{(A_0+H_Q)t} (\lambda- (A_0+H_Q))^{-1}.
\end{equation*}
Since $(\lambda- (A_0+H_Q))e^{(A_0+H_Q)t}$ is a bounded operator on $U$, and $(\lambda- (A_0+H_Q))^{-1}$ is compact for some $\lambda$ by Lemma \ref{l:C0}, it follows that $e^{(A_0+H_Q)t}$ is compact.
\end{proof}

Let $V_0$ be a Hilbert space. An operator $F$ on $V_0$ is said to be diagonalizable if and only if there exists a complete orthonormal system $\{f_n\}_{n \in \mathbb{N}}$ of $V_0$ and a sequence of real numbers $\{\Sigma_n\}_{n \in \mathbb{N}}$ such that
\begin{equation*}
F x= \sum_{n \in \mathbb{N}} \Sigma_n \langle x, f_k \rangle_{V_0} f_n,\quad \forall x \in U.
\end{equation*}
Let $\mathcal{I}$ be an index set. A family of operators $\{F_i\}_{i \in \mathcal{I}} \subseteq V_0$ is said to be diagonalizable simultaneously if and only if there exists a complete orthonormal system $\{f_n\}_{n \in \mathbb{N}}$ of $V_0$ and a sequence of real numbers $\{\Sigma_n^{(i)}\}_{n \in \mathbb{N},i \in \mathcal{I}}$ such that
\begin{equation*}
F_i x= \sum_{n \in \mathbb{N}} \Sigma_n^{(i)} \langle x, f_k \rangle_{V_0} f_n,\quad \forall x \in U , i \in \mathcal{I}.
\end{equation*}

\begin{Lemma}\label{l:diag:ch}
Let $Q$ be such that $\Lambda_{\g,Q}<\infty$ for some $\gamma \in [0, 1/2)$. If $\{E_{k}\}_{k \in \mathbb{N}}$ is the spectral family of $e^{(A_0+H_Q)t_0}$ for some $t_0>0$, then there exists a real-valued sequence  $\{\sigma_k\}_{k \in \mathbb{N}}$ such that for all $t>0$, 
\begin{equation}\label{eq:exp:Sigma}
e^{(A_0+{H}_Q)t}= \sum_{k \in \mathbb{N}} e^{\sigma_k t} E_k. 
\end{equation}
\end{Lemma}
 \begin{proof}
By Lemma \ref{l:analytic},  $e^{(A_0+H_Q)t}$ can be diagonalized  for every $t>0$. We first show that $\{e^{(A_0+H_Q)t}\}_{t \geq 0}$ can be diagonalized simultaneously.
 
 Let $V_k(t)$ be an eigenspace of $e^{(A_0+H_Q)t}$ with eigenvalue $\Sigma_k(t)$. Since the eigenspaces of the compact operator $e^{(A_0+H_Q)t}$ are finite dimensional, $V_k(t)$ is of finite dimension.
 
Let $t_0>0$ be fixed. For any $v \in V_k(t_0)$, we have
\[
e^{(A_0+H_Q)t_0}e^{(A_0+H_Q)t}v = e^{(A_0+H_Q)t}e^{(A_0+H_Q)t_0} v=\Sigma_k(t_0)e^{(A_0+H_Q)t} v,
\]
which implies that $e^{(A_0+H_Q)t}v \in V_k(t_0)$. Therefore, we can restrict $e^{(A_0+H_Q)t}$ to $V_k(t_0)$ for all $t>0$ and obtain a semigroup denoted by $\{e^{(A_0+H_Q)t}|_{V_k(t_0)}\}_{t>0}$ on $V_k(t_0)$. 

Since $V_k(t_0)$ is finite dimensional, $\mathcal{L}(V_k(t_0))$, the space of bounded linear operators on $V_k(t_0)$, is also of finite dimension. By \cite[Corollary 1.3.30]{HC}, $\{e^{(A_0+H_Q)t}|_{V_k(t_0)}\}_{t>0}$ can be diagonalized simultaneously.

Due to the fact that $e^{(A_0+H_Q)t_0}$ is compact and self-adjoint, it holds that $U= \oplus_{k \in \mathbb{N}} V_k(t_0)$, where $\{V_k(t_0)\}_{k \in \mathbb{N}}$ are all the eigenspaces of $e^{(A_0+H_Q)t_0}$. Therefore there exists a complete orthonormal system $\{f_k\}_{k \in \mathbb{N}}$ in $U$ such that
\begin{equation}\label{os}
e^{(A_0+{H}_Q)t} x= \sum_{k \in \mathbb{N}} \Sigma_k(t) \langle x, f_k \rangle_{U} f_k,\quad \forall x \in U, t >0.
\end{equation}

To characterize $ {\Sigma}_{k}(t)$, we observe that ${\Sigma}_{k}(t)$ is differentiable due to the analyticity of the semigroup. Moreover, it satisfies
\begin{equation}\label{func:eq}
{\Sigma}_{k}(t+s) ={\Sigma}_{k}(t){\Sigma}_{k}(s),\quad t, s\geq 0.
\end{equation}
Setting $t=s=0$, we obtain $\Sigma_k(0)=0$ or $\Sigma_k(0)=1$. Since by definition, $e^{(A_0+H_Q)0}$ is the identity operator, it follows that $\Sigma_k(0)=1$.

 By \eqref{func:eq}, we have
 \[
 \frac{{\Sigma}_{k}(t+s) -\Sigma_k(t)}{s} ={\Sigma}_{k}(t)\frac{{\Sigma}_{k}(s)-1}{s},\quad \forall t, s>0
 \]
 Taking the limit as $s\to 0$, we deduce that
 \[
 \Sigma_k(t)'= \Sigma_k(t) \Sigma_k(0)'
 \]
  Hence, there exists $\sigma_{k}\in \mathbb{R}$ such that
\[
{\Sigma}_{k}(t) = e^{{\sigma}_{k} t},\quad \forall k \in \mathbb{N}.
\]

For some fixed $t_0>0$,  the spectral decomposition for $e^{(A_0+{H}_Q)t_0}$ reads 
\[
e^{(A_0+{H}_Q)t_0}= \sum_{k \in \mathbb{N}} e^{\sigma_k t_0} E_k(t_0).
\]
Notice that $e^{\sigma_i t} \neq e^{\sigma_j t}$ for all $t>0$ if and only if $e^{\sigma_i t_0} \neq e^{\sigma_j t_0}$. This observation implies that $\{E_k(t_0)\}_{k \in \mathbb{N}}$ is also the spectral family for any $e^{(A_0+{H}_Q)t}$, $t\geq 0$. Due to the uniqueness of the spectral decomposition, $E_k(t_0)$ does not depend on $t_0$ and \eqref{eq:exp:Sigma} holds.
\end{proof}

Finally, we obtain the following result.

\begin{Corollary}\label{c:diag}
Let $Q$  be such that $\Lambda_{\g,Q}<\infty$ for some $\gamma \in [0, 1/2)$. Let $\{\sigma_{k}\}_{k \in \mathbb{N}}$ and $\{E_k\}_{k \in \mathbb{N}}$ denote the sequence and spectral family introduced in Lemma \ref{l:diag:ch}. Then, the following holds:
\begin{equation}\label{eq:A+H}
A_0 + H_Q= \sum_{k \in \mathbb{N}} \sigma_k E_k.
\end{equation}
\end{Corollary}

\begin{proof}
Let $\{f_k\}_{k \in \mathbb{N}}$ be the complete orthonormal system introduced in \eqref{os}, then we have
\begin{equation*}
e^{(A_0+H_Q)t} f_k = e^{\sigma_k t} f_k.
\end{equation*}
It follows from Lemma \ref{l:analytic} that 
\begin{equation*}
\frac{d}{dt}e^{(A_0+H_Q)t} f_k= (A_0+H_Q)e^{(A_0+H_Q)t} f_k= e^{\sigma_k t}(A_0+H_Q) f_k,
\end{equation*}
and
\begin{equation*}
\frac{d}{dt}e^{(A_0+H_Q)t} f_k= \sigma_k e^{\sigma_k t} f_k.
\end{equation*}
Since $e^{\sigma_k t}>0$, we obtain 
\[
(A_0+H_Q) f_k= \sigma_k f_k, \quad \forall k \in\mathbb{N},
\]
which implies \eqref{eq:A+H}.
\end{proof}

\section{Inverse problem}\label{sec:inverse}

In this section, we address the uniqueness of the inverse problem, which is to determine the covariance operator $Q$ using the correlation data $\theta^{i,j}_{}(t)$ for all $i,j\in\mathbb N$ at a fixed time $t_0$.

\subsection{Cauchy problem for $\theta$}

\begin{Lemma}\label{l:u:regularity}
 Let $Q$ be such that there exists $\g  \in [0, \frac{1}{2}\wedge \frac{d}{4})$ satisfying $ \Lambda_{\g,Q}< \infty$. Then, the following holds:
 \[
 \theta \in L^{\infty}(0,T;U) \cap L^2(0,T;U^1).
 \] 
 Moreover, the correlation function $\theta$ satisfies the following integral equation in $U$:
\begin{equation}\label{eq:U:mild}
 \t(t) = \t(0) + \int_0^t (A_0+ {H}_Q) \t(s) ds,\quad  t \in [0,T].
\end{equation}
\end{Lemma}

\begin{proof}
By Proposition \ref{p:HH}, we obtain 
\[
\|\t(t)\|_U^2= \sum_{i,j \in \mathbb{N}}\E[ \langle u(t) \otimes u(t), e_i\otimes e_j\rangle_U]^2 = \sum_{i,j \in \mathbb{N}} \E[\langle u(t),e_i \rangle \langle u(t), e_j \rangle]^2.
\]
We have from H\"{o}lder's inequality that 
\begin{equation}\label{eq:theta:bounded}
\|\t(t)\|_U^2 \leq \sum_{i,j \in \mathbb{N}} \E[\langle u(t),e_i\rangle^2]\E[\langle u(t),e_j\rangle^2]= \E[\|u(t)\|^2]^2.
\end{equation}
From \eqref{eq:u:reg}, it follows that $\theta \in L^{\infty}(0,T;U)$. A similar argument shows that $\theta \in L^{2}(0,T;U^1)$.

Applying It\^{o}'s formula to $\langle u(t),e_i  \rangle\langle u(t),e_j \rangle$ yields 
\begin{align*}
&\E[\langle u(t),e_i  \rangle\langle u(t),e_j \rangle] = \E[\langle u(0),e_i  \rangle\langle u(0),e_j \rangle] + \E[\int_0^t  \langle A u(s),e_i  \rangle\langle  u(s),e_j \rangle ds] \\
&\quad + \E[\int_0^t\langle  u(s),e_i  \rangle\langle A u(s),e_j \rangle ds] +\E\big[ \int_0^t  \sum_{k \in \mathbb{N}} \langle G(u(s)) Q e_k, e_i \rangle\langle G(u(s)) Q e_k, e_j \rangle ds],\\
\end{align*}
where
\begin{align*}
  \sum_{k \in \mathbb{N}} \langle G(u(s)) Q e_k, e_i \rangle\langle G(u(s)) Q e_k, e_j \rangle &= \sum_{k \in \mathbb{N}} \lambda_k^2 \langle e_k u(s),e_i \rangle \langle e_k u(s),e_j \rangle\\
&= \langle H_Qu(s)\otimes u(s) ,e_i\otimes e_j\rangle_{U}.
\end{align*}
Therefore, we obtain
\begin{align*}
\E[\langle u(t),e_i  \rangle\langle u(t),e_j \rangle]& = \E\left[\langle u(0),e_i  \rangle\langle u(0),e_j \rangle\right] + \E\big[ \int_0^t \langle A_0\big( u(s)\otimes u(s)\big),e_i\otimes e_j  \rangle_U ds \big]  \\
&\quad + \E\big[ \int_0^t  \langle H_Q\big( u(s)\otimes u(s)\big)  ,e_i\otimes e_j\rangle_{U} ds \big]. 
\end{align*}

Similar to the proof of Lemma \ref{l:HQ}, we can deduce that 
\[
\begin{aligned}
&\int_0^t \int_{D \times D} \sum_{k \in \mathbb{N}} \lambda_k^2 |(e_k \otimes e_k) (u(s) \otimes u(s))| |e_i \otimes e_j| dxdy ds\\
&\leq\int_0^t  \Big\|\sum_{k \in \mathbb{N}} \lambda_k^2 |(e_k \otimes e_k)( u(s) \otimes u(s))|\Big\|_{U^{-4\g}} \|e_i \otimes e_j\|_{U^{4\g}}  ds\\
& \leq  \int_0^t  \sum_{k \in \mathbb{N}} \lambda_k^2 \|e_k\|_{d/(2\gamma)}^2 \||u(s)| \otimes  |u(s)|\|_{U} \|e_i \otimes e_j\|_{U^{4\gamma}} ds.
\end{aligned}
\]
By \eqref{eq:theta:bounded}, we have
\[
\begin{aligned}
& \E\Big[\int_0^t \int_{D \times D} \sum_{k \in \mathbb{N}} \lambda_k^2 |e_k \otimes e_k u(s) \otimes u(s)| |e_i \otimes e_j| dxdy ds\Big]\\
& \leq \int_0^t  \Lambda_{\gamma,Q} \E[\| u(s)\|^2] \|e_i \otimes e_j\|_{U^{4\gamma}} < \infty.
\end{aligned}
\]
Thus, applying Fubini's theorem leads to 
\begin{equation*}
\E\big[\int_0^t\langle H_Q\big(u(s)\otimes u(s)\big) ,e_i\otimes e_j\rangle_{U} ds\big]= \int_0^t \langle H_Q \E[u(s) \otimes u(s)], e_i \otimes e_j \rangle_U ds.
\end{equation*}

Similarly, we can show that
\begin{equation*}
 \E[\int_0^t \langle A_0 \big(u(s)\otimes u(s)\big),e_i\otimes e_j  \rangle_U ds] = \int_0^t \langle A_0 \E[ u(s)\otimes u(s)] ,e_i\otimes e_j  \rangle_U ds.
\end{equation*}
Combining the above, we obtain 
\begin{equation*}
\langle \t(t),e_i \otimes e_j \rangle_U  = \langle \t(0),e_i \otimes e_j \rangle_U + \int_0^t \langle (A_0+ {H}_Q) \t(s) ,e_i \otimes e_j \rangle_U ds,\quad  t \in [0,T].
\end{equation*}
Since $\t(t) \in U$ for all $t \in [0,T]$, this implies the integral equation \eqref{eq:U:mild}. 
\end{proof}

\begin{Lemma}\label{l:theta:exp}
 Let $Q$ be such that there exists $\g  \in [0, \frac{1}{2}\wedge\frac{d}{4})$ satisfying $ \Lambda_{\g,Q}< \infty$. Then, for all $\t(0)\in U$, the following holds:
\begin{equation*}
{\t}(t)= e^{(A_0+H_Q)t}\t(0),\quad \forall t\geq 0. 
\end{equation*}
\end{Lemma}

\begin{proof}
 Let $\d(t)= \t(t)-  e^{(A_0+H_Q)t}\t(0)$. Since $e^{(A_0+H_Q)t}\t(0)$ is analytic, by \eqref{eq:U:mild}, we have 
\[
\d(t) = \int_0^t (A_0+H_Q) \d(s) ds.
\]
By Lemma \ref{l:u:regularity}, it follows that $\theta(t) \in L^{\infty}(0,T;U)$. Moreover, since $e^{(A_0+H_Q)t}$ is a strongly continuous semigroup on $U$, we conclude that  $\d(t) \in L^{\infty}(0,T;U)$. 

Recall that $\{S_0(t)\}_{t\geq 0}$ is the semigroup generated by the Laplacian $A_0$. It's well-known \cite[Chapter 5]{DP} that \eqref{eq:d:integral} is equivalent to the following mild form:
\begin{equation}\label{eq:d:integral}
\d(t)= \int_0^t S_0(t-s) H_Q \d(s) ds. 
\end{equation}
By the smoothing effect of the heat kernel \eqref{eq:Schauder} and Lemma \ref{l:HQ}, there holds
\[
\begin{aligned}
\|\delta(t)\|_U=& \int_0^t \|S_0(t-s)\|_{\mathcal{L}(U^{-4\g};U)} \|H_Q\d(s)\|_{U^{-4\gamma}} ds \\
\lesssim& \int_0^t (t-s)^{-2\g} \|H_Q\|_{\mathcal{L}(U;U^{-4\g})} \|\delta(s)\|_U ds.
\end{aligned}
\]

Taking the supreme over $t$ in $[0,t_0]$ gives 
\[
\sup_{t \in[0,t_0]}\|\delta(t)\|_U \lesssim t_0^{1-2\gamma} \Lambda_{\gamma,Q} \sup_{t \in[0,t_0]}\|\delta(t)\|_U,
\]
which implies, for sufficiently small $t_0$, $\sup_{t \in[0,t_0]}\|\delta(t)\|_U=0$, and thus $\delta(t)\equiv 0$ for all $t\in[0, t_0]$. By continuity, we can extend this result to all $t\in [0, T]$, concluding that $\delta(t)=0$ to all $t \in [0,T]$.
\end{proof}

\subsection{Proof of Theorem \ref{thm:main}}

By Lemma \ref{l:theta:exp}, we have
\begin{equation*}%\label{eq:main:1}
\t(t_0,u_0)= e^{(A_0+{H}_Q)t_0}\t(0,u_0).
\end{equation*}
First, we observe that the set $\{\theta^{i,j}(0)\}_{i,j \in \mathbb{N}}$ spans the linear space $U$. Indeed, it can be verified that
\begin{equation*}%\label{eq:main:2}
\begin{aligned}
\langle \theta^{i,j}(0),e_m \otimes e_n \rangle_U &= \langle (e_i+e_j)\otimes (e_i+e_j) - (e_i-e_j)\otimes (e_i-e_j)  ,e_n \otimes e_m \rangle_U\\
&=2 \langle e_j\otimes e_i + e_i\otimes e_j  ,e_n \otimes e_m \rangle_U = 2(\d_{m,i}\d_{n,j}+\d_{m,j}\d_{n,i}).
\end{aligned}
\end{equation*}
Thus, the operator $e^{(A_0+{H}_Q)t_0}$ can be uniquely determined by $\{\theta^{i,j}_{}(t_0)\}_{i,j\in\mathbb N}$. 

From Lemma \ref{l:diag:ch}, there exists a spectral family $\{E_k\}_{k \in \mathbb{N}}$ and a sequence $\{\sigma_k\}_{k \in \mathbb{N}}$ such that
\begin{equation*}
e^{(A_0+{H}_Q)t} = \sum_{k \in \mathbb{N}} e^{\sigma_k t} E_k,
\end{equation*}
where we emphasize that both $\{E_k\}_{k \in \mathbb{N}}$ and $\{\sigma_k\}_{k \in \mathbb{N}}$ can be uniquely determined by $e^{(A_0+{H}_Q)t_0}$. Therefore, the operator $H_Q$ can thus be uniquely determined from the data $\{\theta^{i,j}_{}(t_0)\}_{i,j\in\mathbb N}$ by Corollary \ref{c:diag}. Finally, we observe that
\[
\langle H_Q \chi_D, e_k(x)e_k(y) \rangle_U= \Big\langle \sum_{n \in \mathbb{N}} \lambda_n^2 e_n(x)e_n(y), e_k(x)e_k(y) \Big\rangle_U = \lambda_k^2\ ,\quad k \in \mathbb{N},
\]
which implies that the values $\lambda_k^2$ can be uniquely determined from the data $\{\t^{i,j}_{}\}_{i,j \in \mathbb{N}}$.

\section{Conclusion}\label{sec:conclusion}

In this paper, we investigate the initial boundary value problem for the stochastic heat equation driven by space-time Gaussian noise, which is spatially colored and white in time.  In the inverse problem, we demonstrate that the covariance operator associated with the random potential can be uniquely determined from the correlation of the mild solution at a fixed time. As mentioned in the introduction, a potential direction for future research is to extend the methods developed in this work to the inverse problem for the time-dependent wave equation.


\begin{thebibliography}{10}

\bibitem{CK-SIMA01}
B. Canuto and O. Kavian, Determining coefficients in a class of heat equations via boundary measurements, SIAM J. Math. Anal., 32 (2001), 963--986.

\bibitem{CNYY-IP09}
J. Cheng, J. Nakagawa, M. Yamamoto, and T. Yamazaki, Uniqueness in an inverse problem for a one-dimensional fractional diffusion equation, Inverse Problems, 25 (2009), 115002.

\bibitem{DP}
G. Da Prato and J. Zabczyk,  Stochastic Equations in Infinite Dimensions, Cambridge University Press, 2014.

\bibitem{FLW-IPI24}
X. Feng, P. Li, and X. Wang, An inverse potential problem for the stochastic diffusion equation with a multiplicative noise, Inverse Problems and Imaging, 18 (2024), 271--285. 

\bibitem{HC}
R. A. Horn and C. R. Johnson, Matrix Analysis, Cambridge University Press, 2013

\bibitem{I-98}
V. Isakov, Inverse Problems for Partial Differential Equations, Springer-Verlag, New York, 1998

\bibitem{JZ-IP21}
B. Jin and Z. Zhou, An inverse potential problem for subdiffusion: stability and reconstruction, Inverse Problems, 37 (2021), 015006.

\bibitem{KR-IP19}
B. Kaltenbacher and W. Rundell, On an inverse potential problem for a fractional reaction-diffusion equation, Inverse Problems, 35 (2019), 065004.

\bibitem{KOSY-JDE18}
Y. Kian, L. Oksanen, E. Soccorsi, and M. Yamamoto, Global uniqueness in an inverse problem for time fractional diffusion equations, J. Differential Equations, 264 (2018), 1146--1170.

\bibitem{KF}
A. N. Kolmogorov and S. V. Fomin,  Elements of the Theory of Functions and Functional Analysis, vol. 2: Measure, the Lebesgue Integral, and Hilbert Space, 1961

\bibitem{LPS-CMP08}
M. Lassas, L. P\"{a}iv\"{a}rinta, and E. Saksman, Inverse scattering problem for a two dimensional random potential,
Commun. Math. Phys., 279 (2008), 669--703.

\bibitem{LR}
W. Liu and M. R\"{o}ckner, Stochastic Partial Differential Equations: An Introduction, Springer, 2015.

\bibitem{LLW-SIMA22}
J. Li, P. Li, and X. Wang, Inverse elastic scattering for a random potential, SIAM J. Math. Anal., 54 (2022),5126--5159.

\bibitem{LW-SIMA24}
P. Li and X. Wang, Inverse scattering for the biharmonic wave equation with a random potential, SIAM J. Math.
Anal., 56 (2024), 1959–1995.

%\bibitem{AL}
%A. Lunardi, Analytic Semigroups and Optimal Regularity in Parabolic Problems, Birkh\"{a}user Verlag, 1995.

\bibitem{P}
A. Pazy, Semigroups of Linear Operators and Applications to Partial Differential Equations, Springer, New York, 1983.

\bibitem{RS}
M. Reed and B. Simon, Methods of Modern Mathematical Physics, Academic Press, 1980. 

\end{thebibliography}
\end{document}